\documentclass[11pt,a4paper]{amsart}
\usepackage{amsmath,amsthm,amsfonts,graphicx,enumerate,url,calc,color}
\usepackage[lmargin=27mm,rmargin=27mm,tmargin=27mm,bmargin=27mm]{geometry}
\usepackage[numbers,sort&compress]{natbib}
\renewcommand{\baselinestretch}{1.235}
\newcommand{\half}{\ensuremath{\protect\tfrac{1}{2}}}
\newcommand{\floor}[1]{\ensuremath{\protect\lfloor#1\rfloor}}
\newcommand{\eg}{\textbf{\textup{eg}}}
\newcommand{\ang}[1]{\langle{#1}\rangle}

\newcommand{\thmlabel}[1]{\label{thm:#1}}
\newcommand{\thmref}[1]{Theorem~\ref{thm:#1}}
\newcommand{\lemlabel}[1]{\label{lem:#1}}
\newcommand{\lemref}[1]{Lemma~\ref{lem:#1}}
\newcommand{\twolemref}[2]{Lemmas~\ref{lem:#1} and \ref{lem:#2}}
\newcommand{\eqnlabel}[1]{\label{eqn:#1}}
\newcommand{\eqnref}[1]{\eqref{eqn:#1}}

\newcommand{\eqnrefi}[2]{(\ref{eqn:#1}.#2)}

\newcommand{\figlabel}[1]{\label{fig:#1}}
\newcommand{\figref}[1]{Figure~\ref{fig:#1}}
\newcommand{\seclabel}[1]{\label{sec:#1}}
\newcommand{\secref}[1]{Section~\ref{sec:#1}}

\theoremstyle{plain}
\newtheorem{theorem}{Theorem}
\newtheorem{lemma}[theorem]{Lemma}

\begin{document}

\title{Irreducible Triangulations are Small}

\author{Gwena\"el Joret}
\address{\newline  D\'epartement d'Informatique 
\newline Universit\'e Libre de Bruxelles
\newline Brussels, Belgium}
\email{gjoret@ulb.ac.be}

\author{David~R.~Wood}
\address{\newline Department of Mathematics and Statistics
\newline The University of Melbourne
\newline Melbourne, Australia}
\email{woodd@unimelb.edu.au}

\thanks{
Gwena\"el Joret is a Postdoctoral Researcher of the Fonds National de la Recherche Scientifique (F.R.S.--FNRS).
David Wood is supported by a QEII Research Fellowship from the Australian Research Council.}
\date{\today}

\begin{abstract}
A triangulation of a surface is \emph{irreducible} if there is no edge whose contraction produces another triangulation of the surface.
We prove that every irreducible triangulation of a surface with Euler genus $g\geq1$ has at most $13g-4$ vertices.
The best previous bound was $171g-72$.\\[2ex]
\end{abstract}

\maketitle

\section{Introduction}

Irreducible triangulations are the building blocks of graphs embedded
in surfaces, in the sense that every triangulation can be constructed
from an irreducible triangulation by vertex splitting. Yet there are
only finitely many irreducible triangulations of each surface, as
proved by Barnette and Edelson~\cite{BE-IJM88,BE-IJM89}. Applications of
irreducible triangulations include geometric representations
\cite{BonnNaka-DCG08,ABE-DCG07}, generating triangulations
\cite{Sulanke-Generating,MalnicMohar-JCTB92,NakamotoNegami-DM02,
  LawNeg-JCTB99}, diagonal flips \cite{KNS-JCTB09, Negami-YMJ99,
  CortNaka-DM00}, flexible triangulations \cite{ChenLaw-DM98}, and an
extremal problem regarding cliques in graphs on surfaces \cite{DFJW}. In this paper, we prove
the best known upper bound on the order of an irreducible
triangulation of a surface.

For background on graph theory see \cite{Diestel00}.  We consider
simple, finite, undirected graphs. To \emph{contract} an edge $vw$ in
a graph means to delete $vw$, identify $v$ and $w$, and replace any
parallel edges by a single edge. The inverse operation is called
\emph{vertex splitting}.  Let $G$ be a graph. For a vertex $v\in
V(G)$, let $N_G(v):=\{w\in V(G):vw\in E(G)\}$ and let $G_v$ be the
subgraph of $G$ induced by $\{v\}\cup N_G(v)$.  For $A\subseteq V(G)$,
let $N_G(A):=\bigcup\{N_G(v):v\in A\}$.  Let $e(A)$ be the number of
edges in $G$ with both endpoints in $A$.  For $A,B\subseteq V(G)$, let
$e(A,B)$ be the number of edges in $G$ with one endpoint in $A$ and
one endpoint in $B$.  For $v\in V(G)$, let $e(v,B):=e(\{v\},B)$.

For background on graphs embedded in surfaces see \cite{MoharThom}.
Every surface is homeomorphic to $\mathbb{S}_g$, the orientable
surface with $g$ handles, or to $\mathbb{N}_h$, the non-orientable
surface with $h$ crosscaps. The \emph{Euler genus} of $\mathbb{S}_g$
is $2g$. The \emph{Euler genus} of $\mathbb{N}_h$ is $h$. The
\emph{Euler genus} of a graph $G$, denoted by $\eg(G)$, is the minimum
Euler genus of a surface in which $G$ embeds. A \emph{triangulation}
of a surface $\Sigma$ is a 2-cell embedding of a graph in $\Sigma$,
such that each face is bounded by three edges, and each pair of faces
share at most one edge. A triangulation $G$ of $\Sigma$ is
\emph{irreducible} if there is no edge in $G$ whose contraction
produces another triangulation of $\Sigma$.  Equivalently, for
$\Sigma\neq\mathbb{S}_0$, a triangulation $G$ of $\Sigma$ is
irreducible if and only if every edge of $G$ is in a triangle that
forms a non-contractible cycle in $\Sigma$ \cite{MoharThom}\footnote{
A triangulation of $\Sigma$ is \emph{$k$-minimal} if every
   non-contractible cycle has length at least $k$, and every edge is
   in a non-contractible cycle of length $k$. It is easily seen that
   a triangulation is irreducible if and only if it is
   $3$-minimal. Generalising the result for irreducible
   triangulations, for each surface $\Sigma$ and integer $k$, there
   are finitely many $k$-minimal triangulations of $\Sigma$
   \cite{MN-AMUC95,JMM-JCTB96,GRS-JCTB96}.}.


Recall that there are finitely many irreducible triangulations of each
surface. For example, $K_4$ is the only irreducible triangulation of the sphere
$\mathbb{S}_0$ \cite{SR34}, while $K_6$ and $K_7-E(K_3)$ are the only
irreducible triangulations of the projective plane $\mathbb{N}_1$
\cite{Barnette-JCTB82}. The complete list of irreducible
triangulations has also been computed for the torus $\mathbb{S}_1$
\cite{Lavrenchenko}, the double torus $\mathbb{S}_2$
\cite{Sulanke-Generating}, the Klein bottle $\mathbb{N}_2$
\cite{LawNeg-JCTB97,Sulanke-KleinBottle}, as well as $\mathbb{N}_3$
and $\mathbb{N}_4$ \cite{Sulanke-Generating}.  Gao, Richmond and
Thomassen~\cite{GRT91} proved the first explicit upper bound on the
order of an irreducible triangulation of an arbitrary surface. In
particular, every irreducible triangulation of a surface with Euler
genus $g\geq1$ has at most $(12g+18)^4$ vertices. Nakamoto and
Ota~\cite{NakaOta-JGT95} improved this bound to $171g-72$, which prior
to this paper was the best known upper bound on the order of an
irreducible triangulation of an arbitrary surface. In the case of
orientable surfaces, Cheng~et~al.~\cite{CDP-CGTA04} improved this
bound to $120g$. We prove:

\begin{theorem}
  \thmlabel{Strong} Every irreducible triangulation of a surface with
  Euler genus $g\geq1$ has at most $13g-4$ vertices.
\end{theorem}

The largest known irreducible triangulations of $\mathbb{S}_g$ and of
$\mathbb{N}_h$ respectively have $\floor{\frac{17}{2}g}$ and
$\floor{\frac{11}{2}h}$ vertices \cite{Sulanke06}. Thus the upper
bound in \thmref{Strong} is within a factor of $\frac{26}{11}$ of
optimal.




\section{Background Lemmas}

At the heart of our proof, and that of Nakamoto and
Ota~\cite{NakaOta-JGT95}, is the following lemma independently due to
Archdeacon~\cite{Archdeacon-JGT86} and Miler~\cite{Miller-JCTB87}. Two
graphs are \emph{compatible} if they have at most two vertices in
common.

\begin{lemma}[\cite{Miller-JCTB87,Archdeacon-JGT86}]
  \lemlabel{Miller} If $G$ and $H$ are compatible graphs, then
$$\eg(G\cup H)\geq\eg(G)+\eg(H)\enspace.$$
\end{lemma}

Nakamoto and Ota~\cite{NakaOta-JGT95} proved:

\begin{lemma}[\cite{NakaOta-JGT95}]
  \lemlabel{PositiveGenus} Let $G$ be an irreducible triangulation of
  a surface with positive Euler genus.  Then $G$ has minimum degree at
  least $4$. Moreover, for every vertex $v$ of $G$, the subgraph $G_v$
  has minimum degree at least $4$ and $\eg(G_v)\geq1$.
\end{lemma}

The following definition and lemma is implicit in
\cite{NakaOta-JGT95}. An independent set $S$ of a graph $G$ is
\emph{ordered} if either $S=\emptyset$, or $S$ contains a vertex $v$
such that $S-\{v\}$ is ordered, and $G_v$ and $\bigcup\{G_w:w\in
S-\{v\}\}$ are compatible.  \twolemref{Miller}{PositiveGenus} then
imply:

\begin{lemma}[\cite{NakaOta-JGT95}]
  \lemlabel{Implicit} Let $G$ be an irreducible triangulation of a
  surface with positive Euler genus.  If $S$ is an ordered independent
  set of $G$, then
$$\eg(G)\geq\eg\big(\bigcup_{v\in S}G_v\big)\geq|S|\enspace.$$
\end{lemma}


\section{A Simple Proof}
\seclabel{Simple}

In this section we give a simple proof that every irreducible
triangulation of a surface with Euler genus $g\geq1$ has at most
$25g-12$ vertices. The constant $25$, while greater than the constant
in \thmref{Strong}, is still less than the constant in previous
results. The proof follows the approach developed by Nakamoto and
Ota~\cite{NakaOta-JGT95} (using \lemref{Implicit}). This section also
serves as a helpful introduction to the more complicated proof of
\thmref{Strong} to come.

Let $G$ be an irreducible triangulation of a surface with Euler genus
$g\geq1$.  Let $S$ be a maximal ordered independent set in $G$ such
that $\deg_G(v)\leq6$ for all $v\in S$.  Define
\begin{align*}
  N&:=N_G(S)\text{ ,}\\
  A&:=\{v\in V(G)-(S\cup N):e(v,N)\geq 3\}\text{ ,}\\
  Z&:=\{v\in V(G)-(S\cup N):e(v,N)\leq 2\}\enspace.
\end{align*}
Thus $\{S,N,A,Z\}$ is a partition of $V(G)$.

Suppose that $\deg_G(v)\leq 6$ for some vertex $v\in Z$.  Since
$v\not\in N_G(S)$, the set $S\cup\{v\}$ is independent.  Since
$e(v,N)\leq2$, the subgraphs $G_v$ and $\bigcup\{G_w:w\in S\}$ are
compatible.  Since $S$ is ordered, $S\cup\{v\}$ is ordered.  Hence
$S\cup\{v\}$ contradicts the maximality of $S$.  Now assume that
$\deg_G(v)\geq7$ for all $v\in Z$. Thus
\begin{align}
  \eqnlabel{ZZZ} 7|Z|\leq\sum_{v\in Z}\deg_G(v)= e(N,Z) +
  e(A,Z)+2e(Z)\enspace.
\end{align}
By \lemref{PositiveGenus}, each vertex in $A$ has degree at least $4$,
implying
\begin{align}
  \eqnlabel{AAA} 4|A|\leq\sum_{v\in A}\deg_G(v)= e(N,A) + e(A,Z) +
  2e(A)\enspace.
\end{align}
By the definition of $A$,
\begin{align}
  \eqnlabel{3A} 3|A|\leq\sum_{v\in A}e(v,N)=e(N,A)\enspace.
\end{align}
By Euler's Formula applied to $G$,
\begin{align}
  \eqnlabel{WholeGraph}
  &e(S,N) + e(N) + e(N,A)+ e(N,Z) + e(A) + e(A,Z) + e(Z)\\
  =\; &|E(G)|\;=\; 3(|V(G)|+g-2) \;=\;
  3(|S|+|N|+|A|+|Z|+g-2)\enspace.\nonumber
\end{align}
Summing \eqnref{ZZZ}, \eqnref{AAA}, \eqnref{3A} and
$2\times\eqnref{WholeGraph}$ gives
\begin{align*}
  |A| + |Z|+ 2 e(S,N) + 2e(N) + e(N,Z)\leq 6|S|+6|N|+6(g-2)\enspace.
\end{align*}
Every vertex in $N$ has a neighbour in $S$.  Thus $e(S,N) \geq|N|$.
By \lemref{PositiveGenus}, $G[N]$ has minimum degree at least $3$, and
thus $2e(N)\geq3|N|$.  Since $e(N,Z)\geq0$,
$$|N| +|A| +  |Z|\leq 6|S|+2|N|+6(g-2)\enspace.$$
Since every vertex in $S$ has degree at most $6$, we have
$|N|\leq6|S|$. Thus
$$|V(G)|=|S|+|N|+|A|+|Z|\leq 19|S|+6(g-2)\enspace.$$
By \lemref{Implicit}, $|S|\leq \eg(G)\leq g$.  Therefore $|V(G)|\leq
25g-12$.

\section{Proof of \thmref{Strong}}

The proof of \thmref{Strong} builds on the proof in \secref{Simple}
by:
\begin{itemize}
\item introducing a more powerful approach than \lemref{Implicit} for
  applying \lemref{Miller}, thus enabling \lemref{Miller} to be
  applied to subgraphs with Euler genus possibly greater than $1$
  (whereas \lemref{Implicit} applies \lemref{Miller} to subgraphs with
  Euler genus equal to $1$);
\item choosing an independent set $S$ more carefully than in
  \secref{Simple} so that low-degree vertices are heavily favoured in
  $S$;
\item partitioning $V(G)$ into the similar sets $S,N,A,Z$ as in
  \secref{Simple}, and further partitioning $S$ and $A$ according to
  the vertex degrees;
\item introducing multiple partitions of $N$, one for each value of
  the degree of a vertex in $S$.
\end{itemize}

First we introduce a key definition.  Let $T$ be a binary tree rooted
at a node $r$; that is, every non-leaf node of $T$ has exactly two
child nodes.  Let $L(T)$ be the set of leaves of $T$.  For each node
$x$ of $T$, let $T[x]$ be the subtree of $T$ rooted at $x$.  Suppose
that each leaf $u\in L(T)$ is associated with a given subgraph
$G\ang{u}$ of some graph $G$.  For each non-leaf node $x$ of $T$,
define $$G\ang{x}:=\bigcup_{u\in L(T[x])}\!\!\!G\ang{u}\enspace.$$
Thus $G\ang{x}=G\ang{a}\cup G\ang{b}$, where $a$ and $b$ are the
children of $x$.  The pair $(T,\{G\ang{u}:u\in L(T)\})$ is a
\emph{tree representation} in $G$ if $G\ang{a}$ and $G\ang{b}$ are
compatible for each pair of nodes $a$ and $b$ with a common parent
$x$. In this case, $\eg(G\ang{x})\geq\eg(G\ang{a})+\eg(G\ang{b})$ by
\lemref{Miller}.  This implies the following strengthening of
\lemref{Implicit}.

\begin{lemma}
  \lemlabel{Tree} If $(T,\{G\ang{u}:u\in L(T)\})$ is a tree
  representation in $G$, then
$$\eg(G)\geq \sum_{u\in L(T)}\eg(G\ang{u})\enspace.$$
\end{lemma}

Let $S$ be a set of vertices in a graph $G$.  A tree representation
$(T,\{G\ang{u}:u\in L(T)\})$ in $G$ \emph{respects} $S$ if $L(T)=S$
and $G\ang{u}=G_u$ for each $u\in S$; henceforth denoted
$(T,\{G_u:u\in S\})$.

Let $G$ be an irreducible triangulation of a surface with Euler genus
$g\geq1$.  By \lemref{PositiveGenus}, $G$ has minimum degree at least
$4$.  Let $S$ be an independent set of $G$ such that $\deg_G(v)\leq9$
for all $v\in S$.  For $i\in[4,9]$, define
\begin{align*}
  S_i&:=\{v\in S:\deg_G(v)=i\}\text{ ,}\\
  \widehat{S_i}&:=S_4\cup\dots\cup S_i=\{v\in S:\deg_G(v)\leq i\}\text{ , and}\\
  H_i&:=\bigcup_{v\in\widehat{S_i}}G_v\enspace.
\end{align*}
Observe that $H_4\subseteq H_5\subseteq\dots\subseteq H_9$.  We say
that $S$ is \emph{good} if there is a tree representation
$(T,\{G_u:u\in S\})$ respecting $S$ such that for all $i\in[4,9]$, for
every component $X$ of $H_i$, there is a node $x\in V(T)$ such that
\begin{align*}
  L(T[x])=\widehat{S_i}\cap V(X)\quad\text{ and }\quad
  X=G\ang{x}\enspace.
\end{align*}
Note that these two conditions are equivalent.

For each good independent set $S$ of $G$, let $\phi(S)$ be the vector
$(|S_{4}|, |S_{5}|,\dots,|S_9|)$.  Define
$(a_4,\dots,a_9)\succ(b_4,\dots,b_9)$ if there exists $j \in[4,8]$
such that $a_{i} =b_{i}$ for all $i \in[4,j]$, and $a_{j+1}>b_{j+1}$.
Thus $\succeq$ is a linear ordering.  Hence there is a good
independent set $S$ such that $\phi(S)\succeq\phi(S')$ for every other
good independent set $S'$.  Fix $S$ throughout the remainder of the
proof, and let $(T,\{G_v :v\in S\})$ be a tree representation
respecting $S$.

\begin{lemma}
  \lemlabel{Move} Let $i \in[4,9]$. Suppose that $v$ is a vertex in
  $G-V(H_i)$, such that $\deg_G(v)\leq i$.  Then $v$ has at least
  three neighbours in some component of $H_i$.
\end{lemma}

\begin{proof}
  Suppose on the contrary that $v$ has at most two neighbours in each
  component of $H_i$. Let $j:=\deg_G(v)$.  We now prove that
  $S':=\widehat{S_j}\cup\{v\}$ is a good independent set.

  Say the components of $H_j$ are $X_1,\dots,X_p$, where
  $X_1,\dots,X_q$ are the components of $H_j$ that intersect $N_G(v)$.
  For $\ell\in[1,p]$, the component $X_\ell$ is a subgraph of some
  component of $H_i$.  Thus $v$ has at most two neighbours in
  $X_\ell$.  That is, $G_v$ and $X_\ell$ are compatible.  By
  assumption, for each $\ell\in[1,p]$, there is a node $x_\ell\in
  V(T)$ such that $L(T[x_\ell])=\widehat{S_j}\cap V(X_\ell)$ and
  $X_\ell=G\ang{x_\ell}$.  Thus $T[x_\ell]\cap T[x_k]=\emptyset$ for
  distinct $\ell,k\in[1,p]$.

  Let $T'$ be the tree obtained from the forest
  $\bigcup\{T[x_\ell]:\ell\in[1,p]\}$ by adding a path
  $(v,y_1,\dots,y_p)$, where each $y_\ell$ is adjacent to $x_\ell$.
  Root $T'$ at $y_p$, as illustrated in \figref{Move}.  Observe that
$$L(T')
=\bigcup_{\ell\in[1,p]} \!\!L(T[x_\ell]) \cup\{v\}
=\bigcup_{\ell\in[1,p]} \!\!\big(\widehat{S_j}\cap
V(X_\ell)\big)\cup\{v\} =\widehat{S_j}\cup\{v\}=S'\enspace.$$ Let
$G\ang{u}:=G_u$ for each leaf $u\in L(T')$.  Thus
$G\ang{x_\ell}=X_\ell$ in $T'$, and associated with the node $y_\ell$
is the subgraph $G\ang{y_\ell}=\bigcup\{X_k:k\in[1,\ell]\}\cup G_v$.
The children of $y_1$ are $x_1$ and $v$, and for $\ell\in[2,p]$, the
children of $y_\ell$ are $x_\ell$ and $y_{\ell-1}$.  Since $v$ has at
most two neighbours in $X_\ell$, and since $(X_1\cup\dots\cup
X_{\ell-1})\cap X_\ell=\emptyset$, the subgraphs $G\ang{y_{\ell-1}}$
and $G\ang{x_\ell}$ are compatible.

\begin{figure}
  \begin{center}\includegraphics{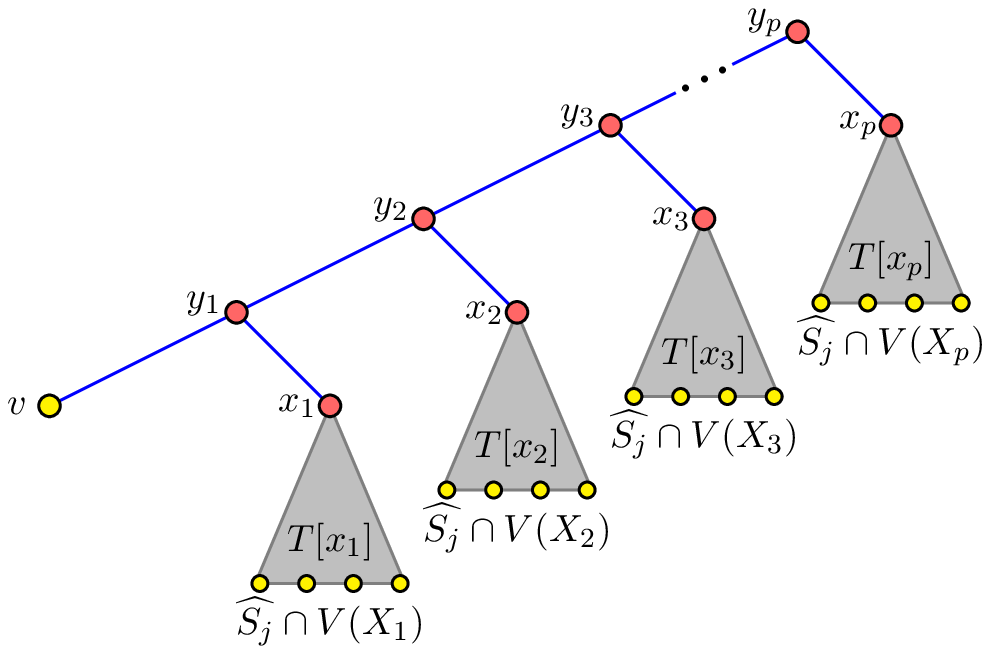}\end{center}
  \caption{\figlabel{Move}Construction of $T'$.}
\end{figure}




Define $H'_4,\dots,H'_9$ with respect to $S'$.  We must prove that for
each $k\in[4,9]$ and for each component $X$ of $H'_k$, there is a node
$z\in V(T')$ such that $X=G\ang{z}$.

First suppose that $k\in[j,9]$.  Since every vertex in $S'$ has degree
at most $j$, we have $\widehat{S_k'} = \widehat{S'_j}$.  Thus $H'_k =
H'_j$, and each component $X$ of $H'_{ k }$ is a component of $H'_j$.
Hence, either $X=X_1\cup \dots\cup X_q\cup G_v$, or $X=X_\ell$ for
some $\ell\in[q+1,p]$.  In the first case, $X=G\ang{y_q}$.  In the
second case, $X=G\ang{x_\ell}$.

Now suppose that $k\in[4,j-1]$.  Thus $S'_k = S_k$ and $H'_k=H_k$.
Each component $X$ of $H'_k$ is a subset of $X_{\ell}$ for some $\ell
\in [1,p]$, and there is a node $z\in T[x_\ell]\subseteq T'$ such that
$X=G\ang{z}$.

This proves that $(T',\{G_u:u\in S'\})$ is a tree representation
respecting $S'$.  Thus $S'$ is a good independent set.  Moreover,
$\phi(S')=(|S_4|,\dots,|S_{j-1}|,|S_j|+1,0,\dots,0)$.  Thus
$\phi(S')\succ\phi(S)$.  This contradiction proves that $v$ has at
least three neighbours in some component of $H_i$.
\end{proof}

\smallskip\noindent\textbf{\boldmath Properties of the Neighbourhood  of $S$:}
Recall that $S$ is a good independent set such that
$\phi(S)\succeq\phi(S')$ for every other good independent set $S'$.
First note that \twolemref{PositiveGenus}{Tree} imply:
\begin{align}
  \eqnlabel{Sg} g\geq \eg(G)\geq \sum_{u\in S}\eg(G_u)\geq|S|\enspace.
\end{align}
Partition $N:= N_G(S)$ as follows. For $i\in[4,9]$, define
\begin{align*}
  U_i&:=N_G(\widehat{S_i})\text{ , }\\
  Y_i&:= N-U_i \text{ , }\\
  V_i &:= \{v \in Y_i: \deg_{G}(v) \leq i\}\text{ , and} \\
  W_i &:= \{v \in Y_i: \deg_{G}(v) \geq i+1\}\enspace.
\end{align*}
Thus $\{U_i,Y_i\}$ and $\{U_i,V_i,W_i\}$ are partitions of $N$ (for
each $i\in[4,9]$).  Also note that $U_4\subseteq
U_5\subseteq\dots\subseteq U_9$, and
$H_i=\bigcup\{G_v:v\in\widehat{S_i}\}$ is a spanning subgraph of
$G[\widehat{S_i}\cup U_i]$.  Each vertex in $N$ has at least one
neighbour in $S$, and each vertex in $S_i$ has $i$ neighbours in
$N$. Thus
\begin{align}
  \eqnlabel{N} |N| \leq e(S,N) = \sum_{i\in[4,9]} i|S_i|\enspace.
\end{align}
For $i\in[5,9]$, each vertex in $U_i-U_{i-1}$ has at least one
neighbour in $S_i$, and each vertex in $S_i$ has at most $i$
neighbours in $U_i-U_{i-1}$. Thus
\begin{align}
  \eqnlabel{Ui} |U_i| \leq |U_{i-1}|+i|S_i|\enspace.
\end{align}
For $i\in[4,9]$, let $c_i$ be the number of components of $H_i$. Thus
$$\sum_{j\in[4,i]}j|S_j|
\geq|E(H_i)| \geq|V(H_i)|-c_i =|U_i|+|\widehat{S_i}|-c_i\enspace.$$
Hence
\begin{align}
  \eqnlabel{Components} |U_i|\leq
  c_i+\sum_{j\in[4,i]}(j-1)|S_j|\enspace.
\end{align}
Consider a vertex $v\in U_i$ for some $i\in[4,9]$.  Thus $v$ is
adjacent to some vertex $w\in\widehat{S_i}$.  It follows from
\lemref{PositiveGenus} that $G[N_G(w)]$ has minimum degree at least 3.
In particular, $v$ has at least three neighbours in $N_G(w)$, which is
a subset of $U_i$.  Thus
\begin{align}
  \eqnlabel{3Ui} 3|U_i| \leq \sum_{v\in U_i} e(v, U_i) = 2e(U_i)
  \enspace.
\end{align}
Consider a vertex $v\in V_i$ for some $i\in[4,9]$. Thus $v$ is in
$G-V(H_i)$ and $\deg_G(v)\leq i$.  By \lemref{Move}, $v$ has at least
three neighbours in $H_i$, implying $e(v,U_i)\geq3$ since
$N_G(v)\cap\widehat{S_i}=\emptyset$.  Hence for $i\in[4,9]$,
\begin{align}
  \eqnlabel{3Vi} 3|V_i| \leq \sum_{v\in V_i} e(v,U_i) = e(U_i,V_i)
  \leq e(U_i,Y_i) \enspace.
\end{align}

\medskip\noindent\textbf{\boldmath Beyond the Neighbourhood of $S$:}
As in \secref{Simple}, partition $V(G)-(S\cup N)$ as
\begin{align*}
  A&:=\{v\in V(G)-(S\cup N):e(v,N)\geq3\}\text{ and}\\
  Z&:=\{v\in V(G)-(S\cup N):e(v,N)\leq2\}\enspace.
\end{align*}
Thus $\{S,N,A, Z\}$ is a partition of $V(G)$.
Further partition $A$ as follows.  For $i\in[4,9]$ let
\begin{align*}
  A_i&:=\{v\in A:\deg_G(v)=i\}\text{, and let}\\
  A_{10}&:=\{v\in A:\deg_G(v)\geq10\}\enspace.
\end{align*}
Thus $\{A_4,\dots,A_{10}\}$ is a partition of $A$.  For $i\in[4,9]$,
let $$\widehat{A_i}:=A_4\cup\dots\cup A_i\enspace.$$

Consider a vertex $v\in A\cup Z$ such that $i=\deg_G(v)\in[4,9]$.  By
\lemref{Move}, $v$ has at least three neighbours in $H_i$, implying
$e(v,N)\geq3$ since $N_G(v)\cap S=\emptyset$.  Thus $v\in A$.  Hence
$\deg_G(v)\geq 10$ for every vertex $v\in Z$.  Since $N_G(Z)\subseteq
A\cup Z\cup N$,
\begin{align}
  \eqnlabel{10Z} 10|Z| \leq \sum_{v\in Z}\deg_G(v) =2e(Z) + e(N,Z) +
  e(A,Z) \enspace.
\end{align}
Since $N_G(A)\subseteq A\cup Z\cup N$,
\begin{align}
  \eqnlabel{iAi} \sum_{i\in[4,10]}i|A_i| \leq \sum_{v\in A}\deg_G(v)
  =2e(A) + e(N,A) + e(A,Z) \enspace.
\end{align}

\medskip\noindent\textbf{\boldmath Global Inequalities:}
Let $i\in[4,9]$.  Consider the sum of the degrees of the vertices in
$Y_i$.  Each vertex in $V_i$ has degree at least $4$, and each vertex
in $W_i$ has degree at least $i+1$.  Each neighbour of a vertex in
$Y_i$ is in $S_{i+1}\cup\dots\cup S_{9}\cup U_{i}\cup Y_{i}\cup A \cup
Z$. Hence
\begin{align*}
  4|V_{i}| + (i+1)|W_{i}| &\leq
  \sum_{v\in Y_i}\deg_G(v)\\
  &\leq
  e(Y_i,S_{i+1}\cup\dots\cup S_{9})+ e(U_{i},Y_{i}) + 2e(Y_{i}) + e(N,Z) + e(Y_i,A)\\
  &\leq \sum_{j\in[i+1,9]}\!\!\! j|S_j| + e(U_{i},Y_{i}) + 2e(Y_{i}) +
  e(N,Z) + e(N,A) -e(U_i,A) \enspace.
\end{align*}
Consider a vertex $v\in\widehat{A_i}$. Thus $v$ is in $G-V(H_i)$ and
$\deg_G(v)\leq i$.  By \lemref{Move}, $v$ has at least three
neighbours in some component of $H_i$, implying $e(v,U_i)\geq3$ since
$N_G(v)\cap S=\emptyset$. Thus
\begin{align*}
  3|\widehat{A_i}| \leq \sum_{v\in\widehat{A_i}}e(v,U_i) \leq e(U_i,A)
  \enspace.
\end{align*}
Hence
\begin{align}
  \eqnlabel{4Vi} 4|V_{i}| + (i+1)|W_{i}| +3|\widehat{A_i}| \leq
  \sum_{j\in[i+1,9]}\!\!\! j|S_j|+ e(U_{i},Y_{i}) + 2e(Y_{i}) + e(N,Z)
  + e(N,A) \enspace.
\end{align}

As proved above, each vertex in $\widehat{A_i}$ has at least three
neighbours in some component of $H_i$.  Let $X_1,\dots,X_{c_i}$ be the
components of $H_i$.  Let $\{D_1,\dots,D_{c_i}\}$ be a partition of
$\widehat{A_i}$ such that for each $\ell\in[1,c_i]$, each vertex in
$D_\ell$ has at least three neighbours in $X_\ell$.  Let $B_\ell$ be
the bipartite subgraph of $G$ with parts $(\widehat{S_i}\cap
V(X_\ell))\cup D_\ell$ and $V(X_\ell)\cap U_i$.

Let $(T,\{G_u:u\in S\})$ be a tree representation respecting $S$.  For
each $\ell\in[1,c_i]$, there is a node $x_\ell$ in $T$ such that
$G\ang{x_\ell}=X_\ell$.  Let $T'$ be the tree obtained from $T$ by
replacing each subtree $T[x_\ell]$ by the single node $x_\ell$.  Thus
$x_\ell$ is a leaf in $T'$.  Redefine $G\ang{x_\ell}:= B_\ell$.  Every
other leaf in $T'$ is a vertex in $S-\widehat{S_i}$.  Thus
$L(T')=(S-\widehat{S_i})\cup\{x_\ell:\ell\in[1,c_i]\}$.  For each
$u\in S-\widehat{S_i}$, leave $G\ang{u}=G_u$ unchanged.  Now
$D_\ell\cap D_k=\emptyset$ for distinct $\ell,k\in[1,c_i]$, and
$G_u\cap D_\ell=\emptyset$ for all $u\in S-\widehat{S_i}$.  Thus
$(T',\{G\ang{u}:u\in L(T')\})$ is a tree representation in $G$.  By
\lemref{Tree},
$$g  \geq \eg(G)\geq 
\sum_{u\in L(T')}\!\!\eg(G\ang{u}) =\sum_{u\in
  S-\widehat{S_i}}\!\!\eg(G_u)\;+\;\sum_{\ell\in[1,c_i]}\!\!\eg(B_\ell)\enspace.$$
By \lemref{PositiveGenus}, $\eg(G_u)\geq1$ for all $u\in
S-\widehat{S_i}$.  Euler's Formula applied to the bipartite graph
$B_\ell$ implies that $|E(B_\ell)|\leq 2(|V(B_\ell)|+\eg(B_\ell)-2)$.
Thus
\begin{align*}
  g &\geq |S-\widehat{S_i}|+\sum_{\ell\in[1,c_i]}\big(\half|E(B_\ell)|
  - |V(B_\ell)| + 2\big)\enspace.
\end{align*}
Since $B_\ell\cap B_k=\emptyset$ for distinct $\ell,k\in[1,c_i]$, and
$\bigcup\{ V(B_\ell):\ell\in[1,c_i]\}=\widehat{S_i}\cup U_i \cup
\widehat{A_i}$,
\begin{align*}
  g &\geq
  |S-\widehat{S_i}|-\big(|\widehat{S_i}|+|U_i|+|\widehat{A_i}|\big)+2c_i+\half\sum_{\ell\in[1,c_i]}|E(B_\ell)|\enspace.
\end{align*}
Each vertex in $\widehat{A_i}$ is incident to at least three edges in
some $B_\ell$.  Thus $$\sum_{\ell\in[1,c_i]}|E(B_\ell)|\geq
3|\widehat{A_i}|+e(\widehat{S_i},U_i)=3|\widehat{A_i}|+\sum_{j\in[4,i]}j|S_j|\enspace.$$
Hence
\begin{align}
  \eqnlabel{Bi}
  g &\geq |S-\widehat{S_i}|-\big(|\widehat{S_i}|+|U_i|+|\widehat{A_i}|\big)+2c_i+\tfrac{3}{2}|\widehat{A_i}|+\half\sum_{j\in[4,i]}j|S_j|\\
  &= |S|-|U_i|+\half|\widehat{A_i}|+2c_i+
  \half\sum_{j\in[4,i]}(j-4)|S_j|\enspace.\nonumber
\end{align}


At this point, the reader is invited to check, using their favourite
linear programming software, that inequalities \eqnref{3A} --
\eqnref{Bi} and the obvious equalities imply that $|V(G)| \leq 13g -
\tfrac{24}{7}$. (Also note that removing any one of these inequalities
leads to a worse bound.)\ What follows is a concise proof of this
inequality, which we include for completeness.

\medskip\noindent\textbf{\boldmath Summing the Inequalities:}
The notation ($x$.$y$) refers to the inequality with label ($x$),
taken with $i=y$.  (For instance, \eqnrefi{Components}{4} stands for
inequality \eqnref{Components} with $i=4$.)

Summing $4 \times $\eqnrefi{Components}{4}, $8 \times
$\eqnrefi{Components}{5}, $4 \times $\eqnrefi{Components}{7}, $4
\times $\eqnrefi{Bi}{4}, $4 \times $\eqnrefi{Bi}{5}, and $2 \times
$\eqnrefi{Bi}{7}, and since $c_{4}\geq 0$,
\begin{align}
  \eqnlabel{step1}
  &4|U_{5}| + 2|U_{7}| + 5|A_{4}| + 3|A_{5}| + |A_{6}| + |A_{7}| + 10|S_{8}| + 10|S_{9}| \\
  \leq\; &38|S_{4}| + 35|S_{5}| + 8|S_{6}| + 11|S_{7}| +
  10g\enspace.\nonumber
\end{align}

Summing $2 \times $\eqnrefi{3Ui}{6}, $4 \times $\eqnrefi{3Ui}{8}, $2
\times $\eqnrefi{3Vi}{6}, $5 \times $\eqnrefi{3Vi}{8}, $2 \times
$\eqnrefi{4Vi}{6}, and $3 \times $\eqnrefi{4Vi}{8} yields
\begin{align*}
  &6|U_{6}| + 12|U_{8}| + 14|V_{6}| + 27|V_{8}| + 14|W_{6}| +
  27|W_{8}| +
  15|A_{4}| + 15|A_{5}| + 15|A_{6}| + 9|A_{7}| + 9|A_{8}|  \\
  \leq\; & 14|S_{7}| + 16|S_{8}| + 45|S_{9}| + 4e(U_{6}) + 8e(U_{8}) +
  4e(U_{6}, Y_{6}) + 8e(U_{8}, Y_{8})
  + 4e(Y_{6}) + 6e(Y_{8}) \\
  &\quad + 5e(N, Z) + 5e(N, A)\enspace.
\end{align*}
Since $6e(Y_{8}) \leq 8e(Y_{8})$ and $e(N) = e(U_{i}) + e(U_{i},
Y_{i}) + e(Y_{i})$ for $i=6$ and $8$, the above inequality becomes
\begin{align}
  \eqnlabel{step2} &6|U_{6}| + 12|U_{8}| + 14|V_{6}| + 27|V_{8}| +
  14|W_{6}| + 27|W_{8}| + 15|A_{4}| + 15|A_{5}| + 15|A_{6}|\\ & +
  9|A_{7}| + 9|A_{8}| \;\leq\; 14|S_{7}| + 16|S_{8}| + 45|S_{9}| +
  12e(N) + 5e(N, Z) + 5e(N, A)\enspace. \nonumber
\end{align}

Summing $12 \times $\eqnrefi{Ui}{6}, $5 \times $\eqnrefi{Ui}{7}, and
$11 \times $\eqnrefi{Ui}{8} gives
\begin{equation}
  \eqnlabel{step3}
  7|U_{6}| + 11|U_{8}| \leq 12|U_{5}| + 6|U_{7}| + 72|S_{6}| + 35|S_{7}| + 88|S_{8}|.
\end{equation} 

Summing \eqnref{step2} and \eqnref{step3}, and since $|N| = |U_{i}| +
|V_{i}| + |W_{i}|$ for $i=6$ and $8$, we obtain
\begin{align*}
  &36|N| + |V_{6}| + 4|V_{8}| + |W_{6}| + 4|W_{8}| + 15|A_{4}| + 15|A_{5}| + 15|A_{6}| + 9|A_{7}| + 9|A_{8}|  \\
  \leq\; & 12|U_{5}| + 6|U_{7}| + 72|S_{6}| + 49|S_{7}| + 104|S_{8}| +
  45|S_{9}| + 12e(N) + 5e(N, Z) + 5e(N, A)\enspace.
\end{align*}
Since trivially $|V_{6}| + 4|V_{8}| + |W_{6}| + 4|W_{8}| \geq 0$ and
$e(N, Z) \geq 0$, the previous inequality implies
\begin{align}
  \eqnlabel{step4}
  &36|N|  + 15|A_{4}| + 15|A_{5}| + 15|A_{6}| + 9|A_{7}| + 9|A_{8}|  \\
  \nonumber \leq\; & 12|U_{5}| + 6|U_{7}| + 72|S_{6}| + 49|S_{7}| +
  104|S_{8}| + 45|S_{9}| + 12e(N) + 6e(N, Z) + 5e(N, A)\enspace.
\end{align}

Summing \eqnref{3A}, $6 \times $\eqnref{10Z} and $6 \times
$\eqnref{iAi} gives
\begin{equation}
  \eqnlabel{step5}
  3|A| + 6\sum_{j \in [4, 10]} j|A_{j}| + 60|Z| \leq 12e(A) + 7e(N, A) + 12e(A, Z) +  12e(Z) + 6e(N, Z)\enspace. 
\end{equation}

Since $|A| = \sum_{i\in [4,10]}|A_{i}|$, summing $3 \times
$\eqnref{step1} with \eqnref{step4} and \eqnref{step5} gives
\begin{align}
  \eqnlabel{step6}
  &57|A| + 3|A_{8}| + 6|A_{10}| + 36|N| + 60|Z|\\
  \nonumber \leq \; &114|S_{4}| + 105|S_{5}| + 96|S_{6}| + 82|S_{7}| + 74|S_{8}| + 15|S_{9}| \\
  \nonumber &\quad + 12e(N) + 12e(A) + 12e(Z) + 12e(N, Z) + 12e(N, A)
  + 12e(A, Z) + 30g\enspace.
\end{align}

Summing $12 \times $\eqnref{WholeGraph} with \eqnref{step6} and since
$e(S, N) = \sum_{i\in[4,9]} i|S_{i}|$, we next obtain
\begin{align*}
  &2|S_{7}| + 22|S_{8}| + 93|S_{9}| + 57|A| + 3|A_{8}| + 6|A_{10}| + 36|N| + 60|Z| \\
  \leq \; & 66|S_{4}| + 45|S_{5}| + 24|S_{6}| + 36|V(G)| + 66g -
  72\enspace.
\end{align*}
Combining this with $4|S_{8}| + 54|S_{9}| + 3|A_{8}| + 6|A_{10}| +
3|Z| \geq 0$, it follows that
\begin{align}
  \eqnlabel{step7}
  &2|S_{7}| + 18|S_{8}| + 39|S_{9}| + 57|A| + 36|N| + 57|Z| \\
  \nonumber \leq \; & 66|S_{4}| + 45|S_{5}| + 24|S_{6}| + 36|V(G)| +
  66g - 72\enspace.
\end{align}

Summing $21 \times $\eqnref{N} with \eqnref{step7} and since $|S| =
\sum_{i\in[4,9]}|S_{i}|$, we derive
\begin{equation*}
  5|S_{7}| + 57|A| + 57|N| +  57|Z| 
  \leq 
  150|S| + 36|V(G)| + 66g - 72.\enspace
\end{equation*}
Since $5|S_{7}| \geq 0$ and $|V(G)| = |S| + |N| + |A| + |Z|$, we
obtain
\begin{equation}
  \eqnlabel{step8}
  21|V(G)|  \leq  207|S| + 66g - 72\enspace.
\end{equation}

Summing $207 \times $\eqnref{Sg} with \eqnref{step8} gives $21|V(G)|
\leq 273g - 72$. That is, $|V(G)| \leq 13g - \frac{24}{7}$, as
claimed. To conclude the proof of \thmref{Strong}, note that
$\frac{24}{7} > 3$, implying that $|V(G)| \leq 13g - 4$ since $|V(G)|$
and $g$ are both integers.


\def\cprime{$'$} \def\soft#1{\leavevmode\setbox0=\hbox{h}\dimen7=\ht0\advance
  \dimen7 by-1ex\relax\if t#1\relax\rlap{\raise.6\dimen7
  \hbox{\kern.3ex\char'47}}#1\relax\else\if T#1\relax
  \rlap{\raise.5\dimen7\hbox{\kern1.3ex\char'47}}#1\relax \else\if
  d#1\relax\rlap{\raise.5\dimen7\hbox{\kern.9ex \char'47}}#1\relax\else\if
  D#1\relax\rlap{\raise.5\dimen7 \hbox{\kern1.4ex\char'47}}#1\relax\else\if
  l#1\relax \rlap{\raise.5\dimen7\hbox{\kern.4ex\char'47}}#1\relax \else\if
  L#1\relax\rlap{\raise.5\dimen7\hbox{\kern.7ex
  \char'47}}#1\relax\else\message{accent \string\soft \space #1 not
  defined!}#1\relax\fi\fi\fi\fi\fi\fi} \def\Dbar{\leavevmode\lower.6ex\hbox to
  0pt{\hskip-.23ex\accent"16\hss}D}

\end{document}